\documentclass[12pt]{amsart}
\usepackage{amsfonts}
\usepackage{latexsym,amscd,amssymb}
\usepackage{movie15}
\usepackage[all]{xy}
\usepackage{array}
\usepackage{hyperref}
\usepackage{graphicx}
\usepackage{epsfig}

\newtheorem*{theorem*}{Theorem}
\newtheorem{theorem}{Theorem}[section]

\newtheorem{lemma}[theorem]{Lemma}

\newtheorem{definition}[theorem]{Definition}
\newtheorem{conjecture}[theorem]{Conjecture}
\newtheorem{example}[theorem]{Example}
\newtheorem*{lemma*}{Lemma}
\newtheorem*{proposition*}{Proposition}
\newtheorem*{corollary*}{Corollary}

\numberwithin{equation}{section}

\allowdisplaybreaks

\begin{document}

\title{Holomorphic Lagrangian fibrations of toric hyperk{\"a}hler manifolds}

\author{Craig van Coevering \and Wei Zhang}

\thanks{Mathematics Classification Primary(2000): Primary 53C26, Secondary 53D20.\\
\indent The second author is supported by Tian Yuan math Fund. and the Fundamental Research Funds for the Central Universities of China\\
\indent Keywords: toric hyperk{\"a}hler manifold, holomorphic
Lagrangian fibration, moment map} \maketitle

\begin{abstract}

For the sake of hyperk{\"a}hler SYZ conjecture, finding holomorphic
Lagrangian fibrations becomes an important issue. Toric
hyperk{\"a}hler manifolds are real dimension $4n$ non-compact
hyperk{\"a}hler manifolds which are quaternion analog of toric
varieties. The $n$ dimensional residue circle action on it
admitting a hyperk{\"a}hler moment map. We use the complex part of this moment map to construct a holomorphic Lagrangian fibration with generic
fiber diffeomorphic to $(\mathbb{C}^*)^n$, and study the singular fibers.

\end{abstract}

\section{Introduction}

In \cite{SYZ96}, Strominger, Yau and Zaslow conjectured that Mirror
Symmetry of Calabi-Yau manifolds comes from real Lagrangian
fibrations. Let $Y$ be a compact, K{\"a}hler, holomorphic symplectic manifold. By
Calabi-Yau theorem, such a manifold admits a hyperk{\"a}hler metric(see \cite{bes87}, or
\cite{huy99}). A complex Lagrangian subvariety of $(Y, I_1)$ is
special Lagrangian with respect to $I_2$, which is clear from the
linear algebra. The hyperk{\"a}hler SYZ conjecture asserts the
existence of holomorphic Lagrangian fibrations on the compact
hyperk{\"a}hler manifolds.

Although the original version is concerning the compact
hyperk{\"a}hler manifold, finding holomorphic Lagrangian fibration
in non-compact hyperk{\"a}hler manifold is also an interesting
problem. For example, Hitchin had constructed holomorphic Lagrangian
fibration in the moduli space of rank-2 stable Higgs bundles of odd
degree with fixed determinant over a Riemann surface(cf.
\cite{hit87self}).

Another important fact must be mentioned is that, in the quest of
examples of special Lagrangian fibration, the physicists set up
Mirror Symmetry in toric Calabi-Yau manifold(cf. \cite{AKMV05},
\cite{Mar10}), and similar topics were also studied by Batyrev(cf.
\cite{bat94}, \cite{bat98}). Thus it is also natural to consider the
fibration in some hyperk{\"a}hler manifold with large symmetry
group.

Toric hyperk{\"a}hler manifolds(some author call them hypertoric
varieties (\cite{Pr08})) are another important class of non-compact
hyperk{\"a}hler manifolds, which are quaternion analogue of toric
varieties. They can be obtained as symplectic quotients of level
sets of the holomorphic moment maps, and themselves admit residue
hyperk{\"a}hler moment maps. Using symplectic quotient technique, in
\cite{BD00}, Bielawski and Dancer studied their moment maps, cores,
cohomologies, etc. While Hausel and Sturmfels study the toric
hyperk{\"a}hler varieties from a more algebraic view
point(\cite{HS02}). Then Konno study them as GIT quotients in \cite{Ko03} and \cite{Ko08}.

We first use the complex
residue moment map to find holomorphic Lagrangian fibration in toric
hyperk{\"a}hler manifold, then study the type of generic and
singular fibers. Namely, let $Y(\alpha,\beta)$ be a toric
hyperk{\"a}hler manifold, then
\begin{theorem*}
The map $\bar{\mu}_\mathbb{C}: Y(\alpha,\beta) \rightarrow
\mathfrak{n}^*_\mathbb{C} \cong \mathbb{C}^n$ defines a holomorphic
Lagrangian fibration, i.e. for any $b \in
\mathfrak{n}^*_\mathbb{C}$, $F_b=\bar{\mu}_\mathbb{C}^{-1}(b)$ is a
complex Lagrangian subvariety.
\end{theorem*}

To study the detail of the generic fiber and singular fiber, for simplicity, we let $\beta=0$. Define a wall structure  $\{W_i\}_{i=1}^d$ on the dual Lie algebra
$\mathfrak{n}^*_\mathbb{C}$(see the precise definition in section
4), there follows

\begin{theorem*}
The generic fiber $F_b$ of $Y(\alpha,\beta)$ over $\mathfrak{n}^*_\mathbb{C} \backslash
\bigcup^d_{i=1} W_i$ is diffeomorphic to complex torus
$(\mathbb{C}^*)^n \cong T^n \times \mathbb{R}^n$.
\end{theorem*}

In the case of $\beta=0$, tt is easy to check that the most singular central fiber $F_0=\bar{\mu}^{-1}_\mathbb{C}(0)$ is the extended core of the toric hyperk{\"a}hler manifold $Y(\alpha,0)$,
constituted by toric varieties intersecting together(cf. \cite{Pr08}). In general, the singular fiber on the discriminant locus is a little complicate, which can be described
by the shrinking torus procedure.

\begin{theorem*}
Consider the singular fiber $F_b$ of $Y(\alpha,\beta)$. When $b$ lies in the generic position of $W_i$, then $F_b$ diffeomorphic to shrinking the real torus $T^1$ generated by $u_i$ in the complex torus $(\mathbb{C}^*)^n \cong \mathbb{R}^n \times T^n$ over the real hyperplane $H_i \subset \mathfrak{n}^* \cong \mathbb{R}^n$. When $b$ lies in the intersection of $s$ walls $\{W_l\}_{l=1}^s$, then the singular fiber is given by shrinking $T^1$ due to $u_l$ over $H_l$ respectively, and at the intersection of $H_{l_i}$, $i=1,\dots,q$, shrinking a torus of real dimension $\dim(\{u_{l_i}\}_{l=1}^q)$ generated by
$\{u_{i_l}\}_{l=1}^q$

\end{theorem*}

The structure of the article is as follows. In section 2, we
introduce some facts of Calabi-Yau and hyperk{\"a}hler geometry, the
special Lagrangian and holomorphic Lagrangian fibration, and background of Mirror Symmetry.

We present the basic properties of toric hyperk{\"a}hler manifold in
section 3. Mainly focus on the symplectic quotient and the GIT
quotient construction. It has close relation with toric variety,
namely the extended core of toric hyperk{\"a}hler manifold are all
constituted by toric varieties and compact toric varieties
respectively, and the cotangent bundle of toric variety in the
extended core is a dense open set of toric hyperk{\"a}hler manifold.

The essential part of this paper is section 4, where we show that
the complex moment map $\bar{\mu}_\mathbb{C}$ defines a holomorphic
Lagrangian fibration. Then we study the type of generic fiber and
singular fiber in great detail.

\

\noindent \textbf{Acknowledgement:} The authors want to thank Prof.
Bin Xu, Prof. Bailin Song and Dr. Yalong Shi for valuable
conversations. The second author is also grateful to prof. Sen Hu
whose string theory lecture in the summer school several years ago
inspired the author's interest in Mirror Symmetry.

\section{Calabi-Yau and hyperk{\"a}hler geometry} A Calabi-Yau
manifold is a K{\"a}hler manifold of complex dimension $n$ with a
covariant constant holomorphic $(n,0)$-form $\Omega$ called the
holomorphic volume form, which satisfies
\begin{equation}
\frac{\omega^n}{n!}=(-1)^{n(n-1)/2}(\sqrt{-1}/2)^n\Omega \wedge
\overline{\Omega},
\end{equation}
where $\omega$ is the K{\"a}hler form. It is a Riemannian
manifold with holonomy contained in $\mathrm{SU}$(n). A (real) Lagrangian subvariety $S$ of an $n$-dimensional Calabi-Yau
manifold is called special Lagrangian if it is calibrated by
$\mathrm{Re}(e^{i\theta}\Omega)$, where $\theta$ is a constant. This condition is equivalent to $\omega|_S=0$ and $\mathrm{Im}
(e^{i\theta}\Omega)|_S=0$(cf. \cite{joy01}).

In 1996, Strominger, Yau and Zaslow \cite{SYZ96} suggested a
geometrical interpretation of Mirror Symmetry between Calabi-Yau
3-folds in terms of dual fibrations by special Lagrangian 3-tori,
now known as the SYZ Conjecture.

A $4n$-dimensional manifold is hyperk{\"a}hler if it possesses a
Riemannian metric $g$ which is K{\"a}hler with respect to three complex
structures $I_1$; $I_2$; $I_3$ satisfying the quaternionic relations
$I_1 I_2 = -I_2 I_1 = I_3$ etc, thus has three forms $\omega_1$,
$\omega_2$, $\omega_3$ corresponding to the three complex
structures. It has holonomy group contained in $\mathrm{Sp}(n)
\subset \mathrm{SU}(2n)$, a prior is Calabi-Yau. With respect to the
complex structure $I_1$ the form $\omega_\mathbb{C} = (\omega_2+
\sqrt{-1} \omega_3)$ is a holomorphic symplectic form. If $X$ is a
complex Lagrangian submanifold ($X$ is a complex submanifold
and $\omega_\mathbb{C}$ vanishes on $X$), then the real and
imaginary parts of $\omega_\mathbb{C}$, namely $\omega_2$ and
$\omega_3$, vanish on $X$. Thus $\omega_2$ vanishes and if $n$ is
odd (resp. even), the real (resp. imaginary) part of $\Omega =
(\omega_3+ \sqrt{-1}\omega_1)^n$ vanishes. Using the complex
structure $I_2$ instead of $I_1$, we see that $X$ is special
Lagrangian(see also \cite{hit97}). A fibration of hyperk{\"a}hler
manifold with complex lagrangian fibers is called holomorphic
Lagrangian fibration, which is also very important in Mirror
Symmetry. This is because of examples of special Lagrangian
fibrations are very rare, all known examples are derived from
holomorphic Lagrangian fibrations on $K3$, torus, or other
hyperk{\"a}hler manifolds. In the weakest form, the hyperk{\"a}hler
SYZ conjecture is stated as follows.

\begin{conjecture}
Let $Y$ be a hyperk{\"a}hler manifold. Then $Y$ can be deformed to a
hyperk{\"a}hler manifold admitting a holomorphic Lagrangian
fibration.
\end{conjecture}
For a more precise form of hyperk{\"a}hler SYZ conjecture, see
\cite{Ver10}.

\section{Toric hyperk{\"a}hler manifold}

One of the most powerful technique for constructing hyperk{\"a}hler
manifolds is the hyperk{\"a}hler quotient method of Hitchin,
Karlhede, Lindstr{\"o}m and Ro{\v c}ek(\cite{HKLR87}). We specialized on the
class of hyperk{\"a}hler quotients of flat quaternionic space
$\mathbb{H}^d$ by subtori of $T^d$. The geometry of these spaces
turns out to be closely connected with the theory of toric
varieties.

Since $\mathbb{H}^d$ can be identified with $T^*\mathbb{C}^d \cong
\mathbb{C}^d\times\mathbb{C}^d$, it has three complex structures
$\{I_1, I_2, I_3\}$. The real torus $T^d=\{(\zeta_1, \zeta_2,
\cdots, \zeta_d) \in \mathbb{C}^d, |\zeta_i|=1\}$ acts on
$\mathbb{C}^d$ induce a action on $T^*\mathbb{C}^d$ keeping the
hyperk{\"a}hler structure,
\begin{equation}
(z,w)\zeta=(z \zeta, w \zeta^{-1}).
\end{equation}

Denote
$M$ the $m$-dimensional connected subtorus of $T^d$ whose Lie
algebra $\mathfrak{m} \subset \mathfrak{t}^d$ is generated by
integer vectors(which is always taken to be primitive), then we have
the following exact sequences
\begin{equation*}
0\rightarrow \mathfrak{m}\xrightarrow{\iota}
\mathfrak{t}^d\xrightarrow{\pi} \mathfrak{n}\rightarrow 0,
\end{equation*}
\begin{equation*}\label{eq:exact:dual}
0\leftarrow {\mathfrak{m}}^*\xleftarrow{\iota^*}
(\mathfrak{t}^d)^*\xleftarrow{\pi^*} {\mathfrak{n}}^*\leftarrow 0,
\end{equation*}
where $\mathfrak{n}=\mathfrak{t}^d/\mathfrak{m}$ is the Lie algebra
of the $n$-dimensional quotient torus $N=T^d/M$ and $m+n=d$.

Let $\{e_i\}_{i=1}^d$ be the standard basis of $\mathfrak{t}^d$ and $\{\theta_i\}_{i=1}^m$ some basis
span $\mathfrak{m}$. Denote $\{e^*_i\}_{i=1}^d$ and $\{\theta^*_i\}_{i=1}^m$ the
dual basis. The subtorus $M$ action admits a hyperk{\"a}hler moment map:
$\mu=(\mu_\mathbb{R}, \mu_\mathbb{C}): \mathbb{H}^d \rightarrow
\mathfrak{m}^* \times \mathfrak{m}^*_\mathbb{C}$, given by,
\begin{equation}\label{eq:moment:r}
\mu_\mathbb{R}(z,w) =
\frac{1}{2}\sum_{i=1}^d(|z_i|^2-|w_i|^2)\iota^* e^*_i,
\end{equation}
\begin{equation}\label{eq:moment:c}
\mu_\mathbb{C}(z,w) = \sum_{i=1}^d z_i w_i \iota^* e^*_i.
\end{equation}
The complex moment map $\mu_\mathbb{C} : \mathbb{H}^d \rightarrow
\mathfrak{m}^*_\mathbb{C}$ is holomorphic with respect to $I_1$.
Bielawski and Dancer introduced the definition of toric
hyperk{\"a}hler varieties, and generally speaking, they are not
toric varieties.

\begin{definition}[\cite{BD00}]
A toric hyperk{\"a}hler variety $Y(\alpha, \beta)$ is a
hyperk{\"a}hler quotient $\mu^{-1}(\alpha, \beta)/M$ for $(\alpha,
\beta) \in \mathfrak{m}^* \times \mathfrak{m}^*_\mathbb{C}$.
\end{definition}

A smooth part of $Y(\alpha,\beta)$ is a $4n$-dimensional
hyperk{\"a}hler manifold, whose hyperk{\"a}hler structure is denoted
by $(g, I_1, I_2, I_3)$. The quotient torus $N=T/M$ acts on
$Y(\alpha,\beta)$, preserving its hyperk{\"a}hler structure. This
residue circle action admits a hyperk{\"a}hler moment map
$\bar{\mu}=(\bar{\mu}_{\mathbb{R}}, \bar{\mu}_{\mathbb{C}})$,
\begin{equation}\label{eq:rmoment:r}
\bar{\mu}_\mathbb{R}([z,w]) =
\frac{1}{2}\sum_{i=1}^d(|z_i|^2-|w_i|^2) e^*_i,
\end{equation}
\begin{equation}\label{eq:rmoment:c}
\bar{\mu}_\mathbb{C}([z,w]) = \sum_{i=1}^d z_i w_i  e^*_i.
\end{equation}
Differs from the toric case, the map $\bar{\mu}$ to $\mathfrak{n}^*
\times \mathfrak{n}^*_\mathbb{C}$ is surjective, never with a
bounded image.

In this article, we always assume that $Y(\alpha,
\beta)$ is a smooth manifold(readers
could consult the regularity argument for \cite{BD00} and
\cite{Ko08}).

\section{\texorpdfstring{$(\mathbb{C^*})^n \cong T^n \times \mathbb{R}^n$}{complex torus}  fibration over \texorpdfstring{$\mathbb{C}^n$}{complex Euclidean space}}

We focus on  the complex moment map $\bar{\mu}_\mathbb{C}$ of the
residue circle action. The first task of this section is to prove
the general theorem.

\begin{theorem}
The map $\bar{\mu}_\mathbb{C}: Y(\alpha,\beta) \rightarrow
\mathfrak{n}^*_\mathbb{C} \cong \mathbb{C}^n$ defines a holomorphic
Lagrangian fibration, i.e. for any $b \in
\mathfrak{n}^*_\mathbb{C}$, $F_b=\bar{\mu}_\mathbb{C}^{-1}(b)$ is a
complex Lagrangian subvariety.
\end{theorem}

\begin{proof}
First of all, $\bar{\mu}_\mathbb{C}([z,w]) = \sum_{i=1}^d z_i w_i
e^*_i$ is obviously a holomorphic map from $Y(\alpha, 0)$ to
$\mathfrak{n}^*_\mathbb{C}$. Let $X$ be the
tangent vector space of $F_b$, for $\bar{\mu}_\mathbb{C}(F_b)=b$, we
have $\mathrm{d} \bar{\mu}_\mathbb{C}(X)=0$. Denote $V$ the $n$
dimensional space of tangent vectors to the orbit of residue circle
action $N$, this is equivalent saying
$\omega_\mathbb{C}(V,X)=\omega_2(V,X)+\sqrt{-1}\omega_3(V,X)=0$,
i.e. $g\langle I_2 V, X\rangle=0$ and $g\langle I_3 V, X\rangle=0$.
At a smooth point of $F_b$, $X$ has real dimension $2n$, thus is the
orthogonal component of $I_2 V \oplus I_3 V$, thus $X$ must be $V
\oplus I_1 V$. By the quaternionic relation, we have $g\langle I_2
X, X\rangle=0$ and $g\langle I_3 X, X\rangle=0$, hence
$\omega_\mathbb{C}|_X=0$, which means that $F_b$ is a complex
Lagrangian subvarieties.
\end{proof}

For the study of the singularity of the fibers, we have to investigate the regularity of the complex moment map
$\bar{\mu}_\mathbb{C}$ of toric hyperk{\"a}hler manifold $Y(\alpha,\beta)$.
We need to reinterpret the dual of Lie algebra $\mathfrak{n}^*$ and $\mathfrak{n}^*_\mathbb{C}$ first. For $\alpha$ lies in $\mathfrak{m}^*$, there is some $x \in \mathbb{R}^d$, such that
\begin{equation}
\alpha=\sum_{i=1}^d x_i \iota^* e_i^*.
\end{equation}
Assume $\iota^* e_i^*=a_i^k \theta_k^*$, and $\alpha=\alpha^k
\theta_k^*$, above equation turns to a linear equation system
\begin{equation}
A x=\alpha,
\end{equation}
where $A$ is $m \times d$ matrix with entry $a_i^k$, and $\alpha$
represents the column vector $\{\alpha^k\}_{k=1}^m$. Its $n$-dimensional solution space is denoted as $\mathfrak{N}_\alpha$, a $n$-plane in
$\mathbb{R}^d$. We can identify $\mathfrak{N}_\alpha$ with $\mathfrak{n}^*$. A hyperplanes arrangement $\{H_i\}_{i=1}^d$ is defined by
$H_i=\mathfrak{N}_\alpha \bigcap \{x_i=0\}$, where
$\{x_i=0\}$ are the coordinate hyperplanes in $\mathbb{R}^d$. Similarly, we could define the $\mathfrak{N}_\beta$ the complex
solution space which identifies to $\mathfrak{n}^*_\mathbb{C}$, and a complex hyperplanes arrangement by    $W_i=\mathfrak{N}_\beta \bigcap \{y_i=0\}$, where $\{y_i=0\}$ is the coordinate hyperplanes in $\mathbb{C}^d$. We call the union $\bigcup^d_{i=1} W_i$ a wall structure on $\mathfrak{N}_\beta$.

\begin{lemma}
The set of regular value of the moment map $\bar{\mu}_\mathbb{C}$ of
toric hyperk{\"a}hler manifold $Y(\alpha, 0)$ is
$$\mathfrak{n}^*_{\mathbb{C}reg}=\mathfrak{n}^*_\mathbb{C} \backslash \bigcup^d_{i=1} W_i.$$
\end{lemma}

\begin{proof}
Let $f: \mathbb{C}^2\rightarrow \mathbb{C}$ be a map defined by
$f(a,b)=ab$. We can easily observe that $(a,b) \in \mathbb{C}^2$ is
a regular point of $f$ if and only if $(a,b) \neq 0$. If one of $(z_i, w_i)$
equals to zero, constrained by equation $\sum_{i=1}^d z_i w_i \iota^* e^*_i=\beta$, the image of
$(\mathrm{d}\bar{\mu}_\mathbb{C})([z,w]) = \sum_{i=1}^d
(\mathrm{d}f)_{(z_i,w_i)} \otimes  e_i^*$ at the point $([z,w]) \in Y(\alpha,\beta)$ can not span the whole $\mathfrak{N}_\mathbb{C}$. Thus the point for some $(z_i,w_i)=0$ is a critical
point of $\bar{\mu}_\mathbb{C}$.
\end{proof}

Immediately, we have
\begin{theorem}
The generic fiber $F_b$ of $Y(\alpha,\beta)$ over $\mathfrak{n}^*_\mathbb{C} \backslash
\bigcup^d_{i=1} W_i$ is diffeomorphic to complex torus
$(\mathbb{C}^*)^n \cong T^n \times \mathbb{R}^n$.
\end{theorem}

\begin{proof}
By the
regularity of the moment map $\bar{\mu}_\mathbb{C}$, $F_b$ is a
smooth manifold. We claim that the residue circle action $N$ acting
on $F$ freely. To see this, lift the $Y(\alpha,\beta)$ to
$\mathbb{H}^d$. Then $[z,w] \in Y(\alpha,\beta)$ has nontrivial isotropy
group in $N$ if and only if the orbit of a subgroup of $N$ through
$(z,w)$ lies in the $M$-orbit. For $N$ is the quotient group,
the only possibility is that some $(z_i, w_i)$ equals to zero, which
can not hold if $b \in \mathfrak{n}^*_\mathbb{C} \backslash
\bigcup^d_{i=1} W_i$. At another hand, the real moment map
$\bar{\mu}_\mathbb{R}$ restricted to $F_b$ is still surjective on
$\mathbb{R}^n$, moreover since $(z_i,w_i) \neq 0$,
$\bar{\mu}_\mathbb{R}$ is also regular, thus $F_b$ must
diffeomorphic to $T^n \times \mathbb{R}^n \cong (\mathbb{C}^*)^n$.
\end{proof}

It is natural to ask what the singular fiber looks like, the
picture will be a little bit complicated. Suggested by the above
proof, we need to investigate the isotropy group in detail.

We first check the simplest case, where $b$ lies in generic
position of $W_i$. Fixing $i$, based on above discussion, the point $[z,w] \in F_b$ where $z_i w_i=0, z_j w_j \neq 0 \
\text{for} \ j \neq i$ but $(z_i,w_i) \neq 0$, is with trivial
isotropy, thus a smooth point on $F_b$. And the real torus $N$ acts
on the complex subvariety $P_{b,i}$ of $F_b$
defined by $(z_i,w_i)=0$ must has a 1 dimensional isotropy subgroup.
For the real moment map restricted to $F_b$ is always surjective,
this is equivalent to shrinking the torus $T^1$ whose Lie algebra is
$u_i$ over the $n-1$ dimensional real subvariety
$\bar{\mu}_\mathbb{R}(P_{b,i})=H_i \subset \mathfrak{n}^* \cong \mathbb{R}^n$, where
$H_i$ is the hyperplane in the real arrangement.

Consider $b$ lies in the intersection of $s$ walls
$\{W_l\}_{l=1}^s$, the situation becomes more
complicated. On the subvariety $P_{b,l}=\{[z,w] \in F_b|(z_l,w_l)=0\}$, we shrink the torus $T^1$ corresponding to $u_l$ over the hyperplane $H_l$. Some of these subvarieties may intersect, or equivalent saying, the hyperplanes will intersect via the real residue moment map. For simplicity, let first $q \leq s$ subvariety intersects. Their normals $\{u_l\}_{l=1}^q$ may be linear dependent, thus we denote $\dim(\{u_l\}_{l=1}^q)$ the dimension of the subspace they spanning. By the Equation (\ref{eq:rmoment:r}), we
know that the image of the real residue moment map is the intersection of $H^l$, $l=1,\dots,q$. Over this intersection, we shrink the
real $\dim(\{u_i\}_{l=1}^q)$ dimensional torus generated by
$\{u_i\}_{l=1}^q$.

Finally, we get the conclusion

\begin{theorem}
Consider the singular fiber $F_b$ of $Y(\alpha,\beta)$. When $b$ lies in the generic position of $W_i$, then  $F_b$ diffeomorphic to shrinking the real torus $T^1$ generated by $u_i$ in the complex torus $(\mathbb{C}^*)^n \cong \mathbb{R}^n \times T^n$ over the real hyperplane $H_i \subset \mathfrak{n}^* \cong \mathbb{R}^n$. When $b$ lies in the intersection of $s$ walls $\{W_l\}_{l=1}^s$, then the singular fiber is given by shrinking $T^1$ due to $u_l$ over $H_l$ respectively, and at the intersection of $H_{l_i}$, $i=1,\dots,q$, shrinking a torus of real dimension $\dim(\{u_{l_i}\}_{l=1}^q)$ generated by
$\{u_{i_l}\}_{l=1}^q$

\end{theorem}

Theoretically, using the data defining the $Y(\alpha,\beta)$, checking
all the intersections of the singular subvarieties, we can
identify the type of the singular fiber. As we know, when $s$ grows
bigger, the computation becomes overwhelming.

\begin{example}
Let the subgroup $M$ generated by $\iota \theta_1=e_1+e_2$ and $\iota \theta_2=e_1+e_3$, then $\iota^* e_1^*=\theta_1^*+\theta_2^*$,
$\iota^* e_2^*=\theta_1^*$, $\iota^* e_3^*=\theta_2^*$. If we take $\alpha=\iota^* e_1- \frac{1}{2}\iota^*e_2=\frac{1}{2}\theta_1^*+\theta_2^*$ and $\beta=0$, consider the toric hyperk{\"a}hler manifold $Y(\alpha,\beta)$ with generic fiber $\mathbb{C}^*$. The
complex equation becomes
\begin{equation*}
\begin{pmatrix}
1 & 1 & 0\\
1 & 0 & 1\\
\end{pmatrix}
\begin{pmatrix}
y_1\\
y_2\\
y_3\\
\end{pmatrix}
=
\begin{pmatrix}
0\\
0\\
0\\
\end{pmatrix},
\end{equation*}
has solution $y=(1,-1,-1)^T$. 

The solution space $\mathfrak{N}_\mathbb{C}$ only intersects the
coordinates hyperplanes in the origin, thus the central fiber is the
only singular fiber. We investigate it closer. For the residue
circle action is $1$-dimensional, the point will be fixed point if
it has nontrivial isotropy. By the defining Equation
(\ref{eq:moment:r}), (\ref{eq:moment:c}) and (\ref{eq:rmoment:c}),
the central fiber $F_0$ satisfies
\begin{equation}
\begin{split}
&|z_1|^2-|w_1|^2+|z_2|^2-|w_2|^2=\frac{1}{2}\\
&|z_1|^2-|w_1|^2+|z_3|^3-|w_3|^2=1\\
&z_1w_1=z_2w_2=z_3w_3=0
\end{split}.
\end{equation}
For 0 is the intersection of 3 walls, consider them respectively: if
$z_1=w_1=0$, then $z_2$ and $w_3$ must be nonzero, similarly, if
$z_2=w_2=0$, then $z_1$ and $w_3$ must be nonzero, if $z_3=w_3=0$,
then $z_2$ and $w_1$ must be nonzero. These are the only 3 fixed
points of $N$. Shrinking a real torus $T^1$ in these 3 points, we
get $\mathbb{C}^1$, $\mathbb{C}P^1$, $\mathbb{C}P^1$ and
$\mathbb{C}^1$ intersecting sequentially, which is nothing but the
toric varieties in the extended core.

\end{example}

Recall that in the category of $T^n$ fibration of toric varieties,
the $T^n$ fibers degenerate at the boundary of the Delzant polytopes
of the toric varieties(cf. \cite{Bou07}). Our theorem can be viewed
as the hyperk{\"a}hler analog. Using moment map to construct special
Lagrangian variety had already been studied by Joyce intensively in
\cite{joy02}.

\bibliographystyle{alpha}
\bibliography{hkBib}

\vfill

\noindent Craig van Coevering: craigvan@ustc.edu.cn

\noindent Wei Zhang: zhangw81@ustc.edu.cn

\noindent School of Mathematics

\noindent University of Science and Technology of China

\noindent Hefei, 230026, P. R.China

\end{document}